\documentclass{article}
\usepackage{amsmath,amssymb,graphicx}
\newcommand{\qed}{\hfill\vrule height6pt  width6pt
depth0pt \medskip}
\title{Explicit formulas for hook walks on continual
Young diagrams
}
\author{}
\date{{\Large Dan Romik} \\
 {\small Laboratoire de Probabilit\'es \\
   Universit\'e Paris 6\\
   e-mail: romik@ccr.jussieu.fr} 
}
\begin{document}
\maketitle
\begin{abstract}
We consider, following the work of S. Kerov, random walks which are
continuous-space generalizations of the Hook Walks defined by
Greene-Nijenhuis-Wilf, performed under the graph of a continual Young
diagram. The limiting point of these walks is a point on the graph of
the diagram. We present several explicit formulas giving the 
probability densities
of these limiting points in terms of the shape of the diagram. This
partially resolves a conjecture of Kerov concerning an explicit formula
for the so-called Markov transform. We also present two inverse formulas,
reconstructing the shape of the diagram in terms of the densities of
the limiting point of the walks. One of the two formulas can be interpreted
as an inverse formula for the Markov transform.
As a corollary, some new integration identities are derived.
\end{abstract}
\linespread{1.9}

\section*{1. Introduction}

A Young diagram is a graphic representation of
a partition $\lambda: \lambda_1 \ge \lambda_2
 \ge ... \ge \lambda_k$ of an integer $n = \sum \lambda_i$. A \emph{continual}
Young diagram is the continous analogue of this, namely a positive
increasing function $f$ defined on some interval $[a,b]$.

Greene, Nijenhuis and Wilf (1979, 1984) introduced two random walks
on Young diagrams
called the \emph{Hook walks}. These random walks continue until reaching
some terminal point on the boundary of the diagram.
By analyzing the probability distributions of these terminal points, 
they reproved 
two important formulas in the combinatorics of Young diagrams. Kerov (1993)
generalized one of the walks to continual Young diagrams. This random walk
converges to a limiting point on the boundary of the diagram, and Kerov
conjectured a formula for the density of this limiting point, in the case
where this random variable is in fact absolutely continuous.

In this paper, we give a unified treatment of both random walks. We prove
Kerov's formula under fairly mild assumptions on the smoothness of the
diagram, and present a new explicit formula for the density of the limiting
point of the second hook walk. The fact that these
expressions are probability densities, and thus integrate to $1$, leads to
some surprising integration relations;
two examples are:
$$ \int_0^1 f(x)dx = $$ $$ \int_0^1 \bigg\{
\frac{1}{\pi} (1+f'(x)) \cdot
\sin\left(\frac{\pi}{1+f'(x)}\right) \cdot $$ $$ \cdot
(x+f(x))^{\frac{f'(x)}{1+f'(x)}} \cdot
(1-x+f(1)-f(x))^{\frac{1}{1+f'(x)}} \cdot $$ $$ \cdot
\exp\left[ -\int_0^1 \frac{1}{u-x+f(u)-f(x)}\cdot
                     \frac{f'(u)-f'(x)}{1+f'(x)} du \right]\bigg\}dx $$
which holds for any positive, increasing, smooth
function $f$ on $[0,1]$ that satisfies $f(0)=0$, and
having first derivative bounded away from
$0$ and infinity and second derivative bounded; and
$$\pi=\int_0^1 \left.{\bigg[} \cos(\pi g(x)/2)\cdot
  x^{-(1+g(x))/2} \cdot (1-x)^{-(1-g(x))/2} \cdot \right.
$$ $$ \left. \exp\left(\frac{1}{2}
        \int_0^1 \frac{g(u)-g(x)}{u-x} du \right){\bigg]}\right.dx
,$$
which holds for any smooth
function $g$ on $[0 ,1]$ which is bounded between $-1+\epsilon$ and 
$1-\epsilon$ for some $\epsilon>0$.

We also solve the inverse problem: that of finding the shape of the
diagram, when given the probability density of the limiting point of
the random walk.  Two inverse formulas are given, one for each of the
two walks. For the walk that was treated by Kerov, the
correspondence between the shape of the diagram and the probability
density of the limiting point of the walk is closely related to the
so-called Markov transform (see Kerov (1998)).
In this case, our explicit formulas
enable the direct calculation of the Markov transform and its inverse.
The Markov transform has found several applications, notably to
Dirichlet priors in statistics, so the inverse formula may well
be applicable to that problem, a possibility that Diaconis and Kemperman
(1996) seem to hint at in their very readable review of the subject.

We remark that the importance of the continual hook walk is best
understood in connection with the asymptotic theory of Plancherel
measure on the symmetric group. Kerov (1999) showed that the
probability density of the limiting point of the walk (the so-called
transition measure - see Section 2 below) governs the
dynamical system of the evolution of a random (Plancherel-distributed)
Young diagram, and used this to illuminate the beautiful 
Vershik-Kerov/Logan-Shepp
limit shape theorem for irreducible representations of the symmetric
group. The transition measure has also appeared in recent work of
Ivanov and Olshanski (2001), where it was shown that the transition
measure of a
random Plancherel-distributed Young diagram converges to the
semicircle distribution, and that the deviation from the semicircle
distribution satisfies a central limit theorem. The form of the
limiting Gaussian process in this central limit theorem exhibits a
surprising resemblance to empirical eigenvalue distribution deviations 
appearing in the GUE random matrix model, a phenomenon that is not yet fully
understood.

In section 2 we give the required definitions and terminology of continual
Young diagrams and the hook walks. We concentrate on ``rotated'' Young
diagrams, so that instead of increasing functions we shall be dealing with
1-Lipschitz functions. However, because of the esthetic appeal in working
with increasing functions, we translate some of the
formulas for those functions.
In section 3 we present the main results, namely the
 formulas for
the densities of the terminating point of the hook walks, together with
the associated integration relations, and the inversion formulas.
We also include a simple asymptotic result on the location of the roots
of the polynomial $\frac{d}{dt}(t(t-1)(t-2)...(t-n))$, which in a sense
inspires the computation for the inverse formulas.

In section 4 we review some of the elementary properties of the hook
walks. The approach using moments is emphasized and some of the results
there may be of independent interest, although the main goal is to
prepare for the proofs of the main results, which are given
in sections 5 and 6. Section 7 contains the formulas for increasing
functions and some more curious formulas related to the hook walks. 

\paragraph{Acknowledgements.} This paper was written during my stay at
the \emph{Laboratoire de Probabilit\'es} of Paris 6 University. I would like
to thank the people of the Laboratoire, and in particular Omer Adelman
and Marc Yor, for their kind encouragement and support. Thanks also to
Philippe Biane for helpful discussions and references.

\section*{2. Definitions}

\paragraph{Continual diagrams.}
While a continual Young diagram is most easily described as an increasing
function on an interval, it turns out that for computational purposes, it is
vastly preferable to use a coordinate system whereby the diagram is
rotated clockwise by an angle of $\pi/4$. We thus define a \emph{diagram},
following Kerov (1993, 1999)
as a 1-Lipschitz function $\omega$ on an interval $[a,b]$, such that 
$a+\omega(a)=b-\omega(b)$. We denote $z=z(\omega)=a+\omega(a)$, the
\emph{center} of the diagram. 
%The point $(z,0)$ shall sometimes be referred to as the corner of the diagram. 
This latter condition makes sure that
the graph of $\omega$ hinges on the graph of the function $x\to|x-z|$,
to make for a true rotated diagram (see Figure 1 below).
Note that equivalently, one may think of a diagram as a 1-Lipschitz function
defined on $\mathbb{R}$, such that outside of some interval $[a,b]$ and
for some $z$, the graph of $\omega$ identifies with $x\to|x-z|$.
The \emph{domain} $D_\omega$ is the set 
$\{(x,y):a\le x\le b, |x-z(\omega)|\le y\le \omega(x)\}$. The \emph{dual
domain} $D_\omega'$ is the set
$\{(x,y):a\le x\le b, \omega(x)\le y\le \min(\omega(a)+x-a,\omega(b)+b-x)\}$
(see Figure 1).
The \emph{area} of $\omega$ is $A(\omega)=\int_a^b (\omega(x)-|x-z|)dx$.
(Note: although this is the true area, it 
is \emph{twice} the area as defined by Kerov (1993, 1999))

\medskip\noindent
Denote by ${\cal D}[a,b]$ the set of diagrams on the interval $[a,b]$.

\paragraph{Smooth diagrams.}
We denote by ${\cal S}[a,b]$ the set of diagrams $\omega\in{\cal D}[a,b]$
satisfying the following smoothness conditions: $\omega$ is piecewise
twice-continuously-differentiable, $\omega''$ is bounded,
and for some two constants $-1<c_1<c_2<1$,
the derivative satisfies $c_1<\omega'(x)<c_2$ whereever it is defined.

\paragraph{Hooks and hook walks.}
For a point $(x,y)\in D_\omega$ ($\omega\in{\cal D}[a,b]$), the
(interior) \emph{hook} of $(x,y)$ is the set
$$\{ (x',y')\in D_\omega: (x'\le x \textrm{\ and\ } y'-y=x-x') \textrm{\ or\ }
                       (x'>x \textrm{\ and\ } y'-y=x'-x)\}.$$
(\emph{In words:} The union of the two rays starting at $(x,y)$ and going
diagonally up-left and up-right, respectively, until they intersect the
graph of $\omega$. 
The intersection with the graph can be a segment,
in which case all the segment is included.)

\medskip\noindent
For a point $(x,y)\in D_\omega'$, the (exterior) hook of $(x,y)$
is the set
$$\{ (x',y')\in D_\omega':(x'\le x\textrm{\ and\ } y-y'=x-x') \textrm{\ or\ }
                       (x'>x\textrm{\ and\ } y-y'=x'-x)\}.$$ 
(\emph{In words:} The union of the two rays starting at $(x,y)$ and going
diagonally down-left and down-right, respectively, until they intersect the
graph of $\omega$.)

\medskip\noindent
The two \emph{hook walks}, the main subjects of this paper, are random
walks on the domain (dual domain, respectively) of the diagram, that,
from a given point, change at each step to a point which is chosen at
random (uniformly, by arc length) from the hook of the last point.

\medskip\noindent
The \emph{exterior corner} walk (or simply: exterior walk)
 starts at the exterior corner point
$(b-\omega(a),\omega(a)+\omega(b))$ and moves at each step to a
uniformly chosen point in the (exterior) hook.

\medskip\noindent
The \emph{interior uniform} walk (or: interior walk)
starts at a uniformly chosen point (by surface area) in
$D_\omega$, and moves at each step to a uniformly chosen point in the
(interior) hook.

\begin{figure}[h]
\begin{center}
\setlength{\unitlength}{0.1 pt}
\begin{picture}(1500,1300)(0,0)
%\put(0,0){\framebox(1500,1300)}
\put(100,100){\vector(1,0){1400}}
\put(100,750){\line(1,-1){650}}
\put(750,100){\line(1,1){650}}

\put(300,120){\line(0,-1){40}}
\put(1350,120){\line(0,-1){40}}
\put(270,0){\shortstack{$a$}}
\put(1330,0){\shortstack{$b$}}
\put(730,0){\shortstack{$z$}}
\put(500,570){\shortstack{$\omega(x)$}}

\multiput(300,550)(30,30){20}{\circle*{1}}
\multiput(1350,700)(-30,30){15}{\circle*{1}}

%\put(300,550){\line(1,1){600}}
%\put(300,550){\line(1,1){100}}
%\put(500,750){\line(1,1){100}}
%\put(700,950){\line(1,1){100}}
%\put(1350,700){\line(-1,1){100}}
%\put(1150,900){\line(-1,1){100}}
%\put(1350,700){\line(-1,1){450}}

\qbezier(300,550)(500,700)(700,650)
\qbezier(700,650)(900,575)(1050,600)
\qbezier(1050,600)(1220,625)(1350,700)
\put(900,1150){\circle*{40}}
\put(1050,1000){\circle*{40}}
\put(950,900){\circle*{40}}
\put(810,760){\circle*{40}}
\put(910,660){\circle*{40}}
\thicklines
\put(900,1150){\line(1,-1){150}}
\put(1050,1000){\line(-1,-1){100}}
\put(950,900){\line(-1,-1){140}}
\put(810,760){\line(1,-1){100}}
\end{picture}
\caption{A continual Young diagram and some steps of an exterior hook walk}
\end{center}
\end{figure}
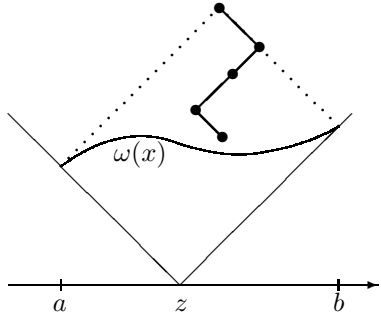

\paragraph{The transition measures.} It is clear that the consecutive steps
of either hook walk must converge almost surely to a limit point
which is on the
graph of the diagram. We call the distribution of the
$x$-coordinate of this limiting point the \emph{transition measure} of the
diagram (relative to the given walk - thus we may talk about the
interior transition measure or exterior transition measure).
The origin of this terminology is in the theory of discrete Young diagrams,
where the transition measures for the exterior-corner- and interior-uniform- walks are in fact the transition measures of some Markov chains, which describe,
respectively, the Plancherel growth of a random Young diagram, and a random
Young tableau of given shape.

\paragraph{Remark.} Note that the starting point of the exterior hook walk
depends on the interval $[a,b]$ where the diagram is defined. This may
cause some confusion in the formulation whereby the diagram is thought of
as a function on $\mathbb{R}$, with the interval $[a,b]$ left unspecified
(the purpose of this formulation was precisely to have a common ground to
discuss diagrams on different intervals, which will be necessary in section
4). However, we remark that the \emph{transition measure} is in fact
\emph{independent} of the choice of interval (as long as the diagram
has its essential support inside the interval, that is, as long as the
diagram identifies with $x\to|x-z|$ outside of the interval). This was
proven for the discrete version of the hook walk in 
Greene-Nijenhuis-Wilf (1979, 1984) - the so-called
``constant zone effect'' - and since, in a sense to be specified in section 4,
the discrete hook walk approximates the general one, the general case
follows. So the choice of interval is in fact immaterial.

\section*{3. The main results}

We now present the main results:

\newpage
\paragraph{Theorem 1. Densities of the transition measures.}
Let $\omega\in{\cal S}[a,b]$. Then:

\medskip \noindent
{\bf (1a)}
 The density of the exterior transition measure for $\omega$ is equal to
$$ 
\frac{1}{\pi} \cos(\pi\omega'(x)/2)\cdot
  (x-a)^{-(1+\omega'(x))/2} \cdot (b-x)^{-(1-\omega'(x))/2} \cdot
 $$ $$ \exp\left(\frac{1}{2}
        \int_a^b \frac{\omega'(u)-\omega'(x)}{u-x} du \right)
$$

\medskip\noindent
{\bf (1b)} The density of the interior transition measure for
 $\omega$ is equal to
$$ \frac{2}{\pi\cdot A(\omega)}
\cos(\pi\omega'(x)/2)\cdot
  (x-a)^{(1+\omega'(x))/2} \cdot (b-x)^{(1-\omega'(x))/2} \cdot
$$ $$ \exp\left(-\frac{1}{2}
        \int_a^b \frac{\omega'(u)-\omega'(x)}{u-x} du \right)
$$

\paragraph{Theorem 2. The continuous ``hook'' integration formulas.}
Let $\omega\in{\cal S}[a,b]$. Then:

\medskip\noindent {\bf (2a)}
$$ \pi = \int_a^b \bigg[
\cos(\pi\omega'(x)/2)\cdot
  (x-a)^{-(1+\omega'(x))/2} \cdot (b-x)^{-(1-\omega'(x))/2} \cdot
 $$ $$ \exp\left(\frac{1}{2}
        \int_a^b \frac{\omega'(u)-\omega'(x)}{u-x} du \right) \bigg] dx $$

\medskip\noindent {\bf (2b)}
$$ \int_a^b (\omega(x)-|x-z(\omega)|)dx  =
$$ $$ = \int_a^b \bigg[ \frac{2}{\pi}
\cos(\pi\omega'(x)/2)\cdot
  (x-a)^{(1+\omega'(x))/2} \cdot (b-x)^{(1-\omega'(x))/2} \cdot
 $$ $$ \exp\left(-\frac{1}{2}
        \int_a^b \frac{\omega'(u)-\omega'(x)}{u-x} du \right) \bigg] dx
$$

\paragraph{Theorem 3. The inversion formulas.} Let $\omega\in{\cal D}[a,b]$.
Then:

\medskip\noindent {\bf (3a)}
If the exterior transition measure of $\omega$ is absolutely continuous,
and its density $g(x)$ is piecewise-continuously-differentiable, has
a bounded derivative, and is bounded away from $0$ (that is,
$\forall x\in[a,b]\ g(x)>c$ for some $c>0$), then for almost all
$x\in[a,b]$
$$ \omega'(x) = -1+\frac{2}{\pi}\textrm{arccot}\left[\frac{1}{\pi}\left(
  \log\left(\frac{b-x}{x-a}\right)
+\frac{1}{g(x)}\int_a^b \frac{g(u)-g(x)}{u-x}du
\right)\right]
$$
(here, and below,
 arccot is the branch of the inverse cotangent function which returns
values between $0$ and $\pi$.)

\medskip\noindent {\bf (3b)}
If the interior transition measure of $\omega$ is absolutely continuous,
and its density $h(x)$ is piecewise-continuously-differentiable, has
a bounded derivative, and is bounded away from $0$, then for almost all
$x\in[a,b]$
$$ \omega'(x) = 1-\frac{2}{\pi}\textrm{arccot}\left[\frac{1}{\pi}\left(
  \log\left(\frac{b-x}{x-a}\right)
+\frac{1}{g(x)}
\left(\int_a^b \frac{g(u)-g(x)}{u-x}du+\frac{2(x-z(\omega))}{A(\omega)}\right)
\right)\right]
$$

As will be seen in section 6, where the proof of Theorem 3 is given,
at the heart of the proof is a limiting calculation involving approximation
of the transition measure by atomic measures.
In the special case where $g(x)$ is the uniform density on $[0,1]$, 
this calculation is particularly simple and may be thought of as a result
on the location of the roots of a certain polynomial. This seems worthy of
mention both for its own sake and as an aid in following the proof of the
general case:

\paragraph{Theorem 4.} Let $p_n(t) = t(t-1)(t-2)(t-3)...(t-n)$. 
The derivative $p_n'$ has a root between each two roots of $p_n$, so
write $p_n'(t) = n\cdot \prod_{k=0}^{n-1} (t-(k+\lambda_{n,k}))$, where
$0<\lambda_{n,k}<1$ are the fractional parts of the roots of $p_n'$. Then we
have for all $0<x<1$,
$$ \lim_{n\to\infty} \lambda_{n,\lfloor x\cdot n\rfloor} =
  \frac{1}{\pi}\textrm{arccot}\left[\frac{1}{\pi}\log\left(\frac{1-x}{x}
\right)\right]
$$
In other words, a plot of the fractional parts of the roots of $p_n'$,
in order of appearance, converges
to a continuous curve. Figure 2 below shows a sample plot of the roots
(in this example, $n=30$) shown against the limiting curve.

\begin{figure}[h]
\begin{center}
\includegraphics{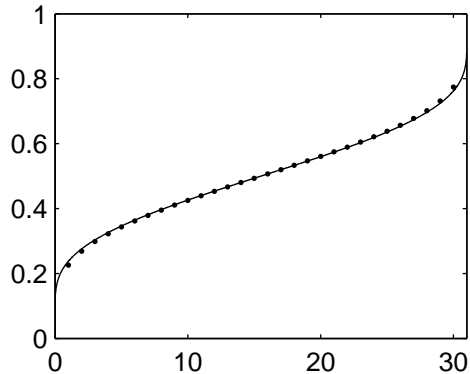}
\caption{Fractional parts of the roots of $p_{30}'$ and the limiting curve}
\end{center}
\end{figure}

\paragraph{Kerov's conjecture.} S. Kerov (1993, 1998) conjectured 
a formula equivalent to formula (1a), as the correct expression not just
for the density of the exterior transition measure in the case when this
measure is absolutely continuous, but more generally for the \emph{
absolutely continuous part} of the exterior transition measure, 
for any diagram in ${\cal D}[a,b]$. Cifarelli and Regazzini (1990)
proved this for \emph{convex} diagrams.
Our Theorem (1a) verifies the conjecture
for the restricted class of diagrams ${\cal S}[a,b]$. However,
we remark that it is quite easy, using the techniques presented in this
paper, to further extend the domain of validity of the formula to
a more general class of diagrams, covering partially the case where the
exterior transition measure is a mixture of an absolutely continuous
part and a discrete part (with no singular component): This is
the class of all the positive, continuous functions
$\omega:[a,b]\to\mathbb{R}$ that are piecewise
twice-continuously-differentiable, that satisfy $a+\omega(a)=b-\omega(b)$,
and such that on any segment of smoothness of $\omega$, either the 
derivative of $\omega$ is bounded between two constants in $(-1,1)$,
or it is identically equal to either $-1$ or $1$. (It is these linear
segments which add the atomic parts to the transition measure.)

\section*{4. Uniqueness, continuity, and moments}

We now review some of the properties of the hook walks on general
diagrams. A special class of diagrams, the \emph{rectangular}
diagrams, will play an important role. These are the diagrams for which
the transition measures are atomic measures with finite support,
so in a sense they are at
the opposite end of the spectrum from the smooth diagrams, and it is using
approximation by these diagrams that the theorems of section 3 will be
proven.

A diagram
$\omega\in{\cal D}[a,b]$ is called rectangular if it is piecewise linear
and its derivative is equal to
$\pm 1$, whereever it exists (Figure 3). Rectangular diagrams have a
particularly simple description using their sets of local minima and maxima:
Let $x_1<x_2<x_3<...<x_n$ be the set of minima of $\omega$, and
$y_1<y_2<...<y_{n-1}$ be its set of maxima. $(x_k)$ and $(y_k)$ are
interlacing sequences, that is, we can write $x_1<y_1<x_2<y_2<...<y_{n-1}<x_n$.
The interlacing sequence pair $(x_k<y_k<x_{k+1})_{k=1}^{n-1}$ determines
$\omega$ uniquely. (In this definition, we think of $\omega$ as a function
on $\mathbb{R}$, identifying outside $[a,b]$ with the function $x\to|x-z|$;
in other words, $a$ and $b$ may not be considered as local maxima, and
may be considered as local minima only if $\omega'(a+)=1$ and
$\omega'(b-)=-1$, respectively.)

The center and the area of a rectangular diagram may be expressed using
the minima and maxima:
$$ z(\omega) = \sum_{k=1}^n x_k - \sum_{k=1}^{n-1} y_k, \qquad
   A(\omega) = 2 \sum_{1\le j\le k\le n-1} (y_j-x_j)(x_{k+1}-y_k) $$
Denote by ${\cal D}_0[a,b]$ the set of all rectangular diagrams on $[a,b]$.

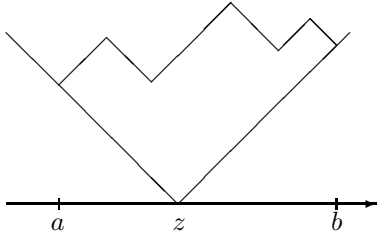
\begin{figure}[h]
\begin{center}
\setlength{\unitlength}{0.1 pt}
\begin{picture}(1500,1000)(0,0)
%\put(0,0){\framebox(1500,1300)}
\put(100,100){\vector(1,0){1400}}
\put(100,750){\line(1,-1){650}}
\put(750,100){\line(1,1){650}}

\put(300,120){\line(0,-1){40}}
\put(1350,120){\line(0,-1){40}}
\put(270,0){\shortstack{$a$}}
\put(1330,0){\shortstack{$b$}}
\put(730,0){\shortstack{$z$}}
%\put(500,570){\shortstack{$\omega(x)$}}

\put(300,550){\line(1,1){180}}
\put(480,730){\line(1,-1){170}}
\put(650,560){\line(1,1){300}}
\put(950,860){\line(1,-1){178}}
\put(1250,800){\line(-1,-1){122}}
\put(1350,700){\line(-1,1){100}}

\end{picture}
\caption{A rectangular diagram}
\end{center}
\end{figure}

\subsection*{4.1. The exterior walk}

The results of this subsection have appeared in Kerov (1993), we include them
for completeness and to motivate the analogous results of the following
subsection, which discusses the interior walk.

Our starting point is the formula for the exterior transition measure of
a rectangular diagram. The transition measure $\mu$
in this case is clearly atomic,
and concentrated on the set of minimum points $(x_k)$. We have for
$k=1,2,...,n$,
$$ (1) \quad
\mu(x_k)
  = \frac{\prod_i (x_k-y_i)}{\prod_{i\neq k}(x_k-x_i)}
  = \prod_{i<k} \left(1-\frac{y_i-x_i}{x_k-x_i}\right) \cdot
    \prod_{i>k} \left(1-\frac{x_i-y_{i-1}}{x_i-x_k}\right)
$$
For the proof see Kerov (1993).
Equivalently, one may define $\mu(x_k)$ using the partial fraction
decomposition
$$ (2) \qquad\qquad \sum_{k=1}^n \frac{\mu(x_k)}{x-x_k} =
\frac{\prod_{i=1}^{n-1} (x-y_i)}{\prod_{i=1}^n (x-x_i)} $$
$(2)$ can be rewritten as
$$ (3) \qquad\qquad \int_\mathbb{R} \frac{d\mu(t)}{x-t} = 
\frac{1}{x}\exp\left(\int_\mathbb{R} \frac{d\sigma(t)}{t-x}\right)$$
where $\sigma$ is the \emph{charge} of the diagram $\omega$,
defined as the function $\sigma(x)=(\omega(x)-|x|)/2$. This holds for real $x
\notin [a,b]$. We now show that $(3)$ can in
fact be taken as an alternative defining equation for the exterior
transition measure of \emph{any} diagram (i.e. not just a rectangular one):

\paragraph{Lemma 1.} For any diagram $\omega\in{\cal D}[a,b]$ with
charge $\sigma$ and exterior transition measure $\mu$, $(3)$ holds.

\paragraph{Proof.} Equip ${\cal D}[a,b]$ with the topology of uniform
convergence, and the set of measures on $[a,b]$ with the weak topology.
Clearly, ${\cal D}_0[a,b]$ is dense in ${\cal D}[a,b]$. Let 
$\omega_n$ be a sequence of rectangular diagrams converging to $\omega$.
It is easy to see that the distribution of the entire exterior
hook walk on $\omega_n$ converges weakly to the distribution of the 
exterior hook walk on $\omega$, and in particular, the transition measures
$\mu_n$ of $\omega_n$ converge to $\mu$. Also, the signed measures $d\sigma_n$
corresponding to the charges of $\omega_n$, converge weakly to $d\sigma$. 
Therefore, for $x\notin[a,b]$
we have $$ \int_\mathbb{R} \frac{d\mu_n(t)}{x-t} 
\xrightarrow[n\to\infty]{} \int_\mathbb{R}\frac{d\mu(t)}{x-t}, $$

$$ \frac{1}{x}\exp\left(\int_\mathbb{R} \frac{d\sigma_n(t)}{t-x}\right)
\xrightarrow[n\to\infty]{}
\frac{1}{x}\exp\left(\int_\mathbb{R} \frac{d\sigma(t)}{t-x}\right) $$
Since $(3)$ holds for each of the $\omega_n$, it holds for $\omega$.
\qed

We now rephrase equation $(3)$ using moments. For a diagram $\omega\in
{\cal D}[a,b]$ with charge $\sigma$ and exterior transition measure $\mu$,
define $p_n=-n\int_\mathbb{R} u^{n-1} d\sigma(u)$
($n=1,2,3,...$), and $h_n = \int_\mathbb{R} u^n d\mu(u)$ ($n=0,1,2,...)$.
By expanding into power series the integrands on both sides of $(3)$,
we can rewrite it as an identity of generating functions
$$ (4)\qquad
\sum_{n=0}^\infty h_n x^{-n} = \exp\left(\sum_{n=1}^\infty \frac{p_n}n
  x^{-n}\right) $$
The two series converge when $|x|$ is large enough. By equating the 
coefficients on both sides one obtains the relation
$$ (5) \qquad h_n = \sum_{\rho \vdash n} \prod_{k\ge 1} 
               \left(\frac{p_k}{k}\right)^{\rho_k}/\rho_k!, $$
where $\rho=(\rho_1,\rho_2,\rho_3,...)$ runs over all partitions of $n$.
($\rho_k$ indicates the number of times $k$ appears in the partition, so
$ n=\sum k\rho_k$.) From this it is easy to see by induction 
that $h_n$ also determines
$p_n$ uniquely (in fact each $p_n$ is a polynomial in $h_1,...,h_n$).

\medskip\noindent
We are now in a position to prove:

\paragraph{Theorem 5.} The correspondence $\omega\to\mu$ which assigns
to a diagram $\omega$ its exterior transition measure $\mu$ establishes
a homeomorphism between the set ${\cal D}[a,b]$ and the set ${\cal M}[a,b]$
of probability
measures on the interval $[a,b]$ (with the topologies defined above).

\paragraph{Proof.} Proofs may be found in Kerov (1993), Krein-Nudelman
(1977). We give a proof which is a variation on Kerov's proof:
If $\omega_n\to\omega$ in ${\cal D}[a,b]$, then, using the same
argument as in the proof of Lemma 1 above, because of weak convergence of
the distribution of the entire hook walk, we also have convergence
$\mu_n\to\mu$ of the transition measures. So the correspondence $\omega
\to\mu$ is continuous. We now prove that it is invertible and its
inverse is continuous: if $\mu\in{\cal M}[a,b]$, take a sequence of
atomic probability measures with finite support $\mu_n\in{\cal M}[a,b]$
converging weakly to $\mu$. For each such $\mu_n$ there exists a 
(unique)
rectangular diagram $\omega_n$ whose transition measure is $\mu_n$
(this follows directly from (2) - the $(x_k)$ are the atoms of $\mu_n$
and the $(y_k)$ are the roots of the equation $\sum \mu(x_k)/(x-x_k) = 0$).
Since the $\omega_n$ are $1$-Lipschitz and satisfy $0\le \omega(a)\le b-a$,
they are equicontinuous and uniformly bounded, therefore by the Arzela-Ascoli
theorem have a uniformly convergent subsequence
${{\omega_n}_k\to\omega}$.
By the continuity proven above, ${\mu_n}_k\to\mu'$ where $\mu'$ is the
transition measure of $\omega$. But ${\mu_n}_k\to\mu$, so $\mu=\mu'$ and
$\omega$ is the desired inverse image of $\mu$. The uniqueness of the
inverse image of $\mu$ follows from the fact stated above that the moments
$h_n$ of $\mu$ determine uniquely the moments $p_n$ of $\sigma$, which
determine $\sigma$ (and therefore $\omega$) uniquely. Finally, if
$\mu_n\to\mu$ and $\omega, \omega_n$ are the inverse images of $\mu, \mu_n$,
respectively, then any subsequence ${\omega_n}_k$ contains (by Arzela-Ascoli)
a convergent
subsequence ${{\omega_n}_k}_j$, which, by uniqueness and continuity, must converge
to $\omega$. Therefore $\omega_n$ itself must converge to $\omega$. This
establishes continuity of the inverse correspondence and finishes the proof.
\qed

\subsection*{4.2. The interior walk}

The interior walk exhibits an interesting duality with the exterior walk,
so the ideas of the previous subsection copy over, with some minor changes,
to the case of the interior walk. One notable complication is that the
correspondence assigning to a diagram its interior transition measure is
not one-to-one, so we do not have uniqueness of an inverse image. However,
uniqueness can be restored if we fix two parameters, the center and the
area of the diagram.

We start, as before, with the formula
for the interior transition measure of a rectangular diagram.
The transition measure $\nu$ will in this case be concentrated on the
maximum points $(y_k)$ of the diagram, with the sizes of the atoms being
$$ (6)\qquad \nu(y_k)
  = -\frac{2}{A(\omega)} \cdot 
 \frac{\prod_i (y_k-x_i)}{\prod_{i\neq k}(y_k-y_i)} = $$ $$
  = \frac{2}{A(\omega)} \cdot (x_{k+1}-y_k)(y_k-x_k)\cdot 
   \prod_{i<k} \left(1+\frac{y_i-x_i}{y_k-y_i}\right) \cdot
   \prod_{j>k} \left(1+\frac{x_{j+1}-y_j}{y_j-y_k}\right) $$
for $k=1,2,...,n-1$.
While formula $(5)$ appears explicitly in Kerov (1993, 1999),
$(6)$ appears in a somewhat different
form in Greene-Nijenhuis-Wilf's (1979) treatment of the discrete hook walk.
To formally deduce it from their result, one needs to consider first
a walk on rectangular diagrams having their center at $0$ and all of
whose local extrema lie on integer points of the plane. For those diagrams,
their formula translates to (6) upon conversion to rotated 
coordinates. Next, it can be seen by scaling reasons
that the formula is valid for diagrams whose local extrema lie
on rational points of the plane. And then, by approximation the result
follows. Alternatively, one may form a Markov chain analogous to the one
in proposition 4.1 of Kerov (1993)
and use arguments similar to the ones in the original paper of
Greene-Nijenhuis-Wilf (1979) (see also Pittel (1986))
to give a direct proof.

As in the exterior walk case, $(6)$ is equivalent to the partial
fraction decomposition
$$ (7)\qquad
-\frac{A(\omega)}{2}\sum_{k=1}^{n-1} \frac{\nu(y_k)}{x-y_k} + x - z(\omega)
= \frac{\prod_{i=1}^n (x-x_i)}{\prod_{i=1}^{n-1} (x-y_i)} $$
(this is one way of verifying that the $\nu(y_k)$ sum to $1$)
which can be rewritten as
$$ (8)\qquad
-\frac{A(\omega)}{2}\int_\mathbb{R} \frac{d\nu(t)}{x-t} + x-z(\omega)=
   x\exp\left(-\int_\mathbb{R}\frac{d\sigma(t)}{t-x}\right) $$
We state for the record:

\paragraph{Lemma 2.} (8) holds for any diagram $\omega\in{\cal D}[a,b]$
with charge $\sigma$ and interior transition measure $\nu$.

\paragraph{Proof.} Take a sequence $\omega_n \in{\cal D}[a,b]$ of
diagrams approximating $\omega$ and having the same area and center as
$\omega$, and continue as in the proof of Lemma 1. \qed

As before, we translate (8) into the language of moments. With
$p_n=-n\int_\mathbb{R} u^{n-1}d\sigma(u)\ \  (n=1,2,...)$ as before, and
$g_n = \int_\mathbb{R} u^n d\nu(u)\ \  (n=0,1,2,...)$, we have the equation
$$ (9) \qquad
 1-\frac{z}{x} - \frac{A}{2}\sum_{n=0}^\infty g_n x^{-(n+2)} =
 \exp\left(-\sum_{n=1}\frac{p_n}{n} x^{-n} \right),
$$
valid for large enough $|x|$. Upon equation of the coefficients one obtains
$$ (10)\quad
g_{n-2}=-\frac{2}{A} \sum_{\rho \vdash n} \prod_{k\ge 1} 
               \left(\frac{-p_k}{k}\right)^{\rho_k}/\rho_k! =
      \frac{2}{A} \sum_{\rho \vdash n} (-1)^{\sum\rho_k+1} \prod_{k\ge 1} 
               \left(\frac{p_k}{k}\right)^{\rho_k}/\rho_k!
$$
(We also get $z=p_1$ and, since $g_0=1$, $A=p_2-p_1^2$.) 
We note as before that
this implies that $\nu$ determines $\omega$ uniquely, \emph{provided the
area and center are fixed}. Another small complication relative to the
case of the exterior walk, is that the support of $\nu$ is generally smaller
than the support of $\omega(x)-|x-z|$.
Also note that the trivial rectangular diagram $x\to|x-z|$ does \emph{not}
have an interior transition measure. This leads us to the following
analogue of Theorem 5:

\paragraph{Theorem 6.} Fix $z\in\mathbb{R}, A> 0$.
Let ${\cal D}_{A,z}[a,b]$ be the set of all diagrams on $[a,b]$ 
having area $A$ and center $z$, equipped with
the topology of uniform convergence. The correspondence $\omega\to\nu$
assigning to a diagram $\omega\in{\cal D}_{A,z}[a,b]$ its interior
transition measure $\nu$ is a homeomorphism between ${\cal D}_{A,z}[a,b]$
and some closed subset of ${\cal M}[a,b]$. In the inverse direction,
for each probability measure $\nu\in{\cal M}[a,b]$, there exists a 
unique diagram
$\omega\in{\cal D}_{A,z}[c,d]$ for some interval $[c,d]\supset[a,b]$ such that
$\nu$ is its interior transition measure.
 $c$ and $d$ may be taken to depend
only on $a$ and $b$, and not on $\nu$.
The correspondence $\nu\to\omega$
is a homeomorphism between ${\cal M}[a,b]$ and some closed subset of 
${\cal D}_{A,z}[c,d]$.

\paragraph{Proof.} This is basically a repetition of the arguments used in
the proof of Theorem 5. The only new fact left to prove, in order to enable
the use of the Arzela-Ascoli theorem and to prove the last claim in Theorem 6,
is the following: For a given area $A$, center $z$, and interval $[a,b]$,
there exists an interval $[c,d]\supset[a,b]$ such that for every atomic
measure $\nu$ with finite support, the unique rectangular diagram $\omega$
having $\nu$ as its interior transition measure (whose existence is guaranteed
by (7)) has its support (or more precisely the support of $\omega(x)-|x-z|$)
contained in $[c,d]$.

We proceed to prove this fact: Let $x_1<y_1<x_2<...<y_{n-1}<x_n$ be
the interlacing sequences of minima and maxima of $\omega$. $y_1,y_2,...,
y_{n-1}$ are the atoms of $\nu$, so it suffices to find an interval $[c,d]$
guaranteed to contain $x_1,x_n$. $x_n$ is the root of the equation
$$ -\frac{A}{2}\sum_{k=1}^{n-1} \frac{\nu(y_k)}{x-y_k}+x-z = 0$$
which lies to the right of $y_{n-1}$. But since
$$ -\frac{A}{2}\sum_{k=1}^{n-1} \frac{\nu(y_k)}{x-y_k}+x-z >
  -\frac{A}{2}\cdot\frac{1}{x-y_{n-1}}+x-z,$$ and both of these functions
are increasing on $(y_{n-1},\infty)$, $x_n$ can be no greater than the
root of the equation
$$ -\frac{A}{2}\cdot\frac{1}{x-y_{n-1}}+x-z = 
\frac{x^2-(y_{n-1}+z)x+(zy_{n-1}-A/2)}{x-y_{n-1}} = 0 $$
Remembering that $y_{n-1}\le b,$ we have $x_n \le (b+z+\sqrt{(b+z)^2+2A})/2$,
which gives us the upper bound $d$ for the support of $\omega$. $d$ depends
only on $a$ and $b$ since we have trivially $a\le z\le b$ and
$0<A\le (b-a)^2/2$, otherwise ${\cal D}_{A,z}[a,b]=\emptyset$. 
The lower bound is obtained in a similar way. \qed

\paragraph{Remark.} Theorem 6 implies that the correspondence 
$\omega\to(A(\omega),z(\omega),\nu)$
gives a bijection between 
${\cal D}_\mathbb{R}:=\cup_{a<b}{\cal D}[a,b]\setminus\{
\textrm{trivial diagrams\ } x\to|x-z|\}$
and $(0,\infty)\times\mathbb{R}\times{\cal M}_\mathbb{R}$,
where ${\cal M}_\mathbb{R}:=\cup_{a<b}{\cal M}[a,b]$. However, it can be
seen that this
bijection is \emph{not} a homeomorphism when ${\cal D}_\mathbb{R}$ and
${\cal M}_\mathbb{R}$ are equipped with the topologies of uniform
and weak convergence, respectively. It is interesting to ask what is
the precise topological nature of this correspondence. For our
purposes, however, the results of Theorem 6 suffice.

\section*{5. Calculation of the densities}

In this section, we calculate the densities of the exterior and interior
transition measures for smooth diagrams. The basic tool is
to approximate smooth diagrams by rectangular ones and use formulas
(1) and (6). We write a detailed analysis of the exterior case,
and go rapidly through the calculation in the interior case.

\subsection*{5.1. The exterior transition measure}

We shall prove Theorem (1a) in two approximation steps. First, we prove it
for diagrams $\omega\in{\cal S}[a,b]$ which are piecewise linear. Then,
we shall approximate arbitrary smooth diagrams by piecewise linear ones. We
shall use the two following well-known relations involving the
gamma function:
$$ (11) \quad
 \prod_{k=1}^n \left(1+\frac{t}{k}\right) \sim \frac{n^t}{\Gamma(t+1)},
\qquad (12) \quad \Gamma(t)\Gamma(1-t) = \frac{\pi}{\sin(\pi t)}$$

\subsubsection*{5.1.1. Piecewise linear diagrams}

Let $\omega\in 
{\cal S}[a,b]$ be piecewise linear. Let $\mu$ be the measure on
$[a,b]$ whose density $g(x)$ is given by the right-hand side of
(1a). We define a sequence
$\omega_n$ of rectangular diagrams approximating $\omega$, as follows:
First note that a rectangular diagram is determined uniquely by giving
its local minima $x_i$ and the values there. Now define $\omega_n$ by
having its local minima be the points dividing each of the segments
of linearity of $\omega$ into $n$ equal parts, together with the
requirement that $\omega_n$ interpolate $\omega$ at these points.

Let $\mu_n$ be the exterior transition measure of $\omega_n$.
Our claim is that $\mu_n \to \mu$ weakly as $n\to\infty$.
We will work within each segment of linearity of $\omega$, and finally
``glue'' the results together.

Let $[A,B]\subset [a,b]$ be a segment of linearity of $\omega$.
Thinking of $n$ as fixed for the moment, let $A=x_N<x_{N+1}<
x_{N+2}< ... < x_{N+n}=B$ be the minima of $\omega_n$ within the segment
$[A,B]$, that is, $x_{N+k} = A+(B-A)k/n, k=0,1,2,...,n$. It is easy
to calculate that $y_{N+k} = A+(B-A)(k+1/2)/n + (B-A)d/2n$, where
$d=(\omega(B)-\omega(A))/(B-A)$ is the slope of $\omega$ on $[A,B]$
($y_{N+k}$ is calculated by reasoning that it
lies on the intersection of the two lines whose equations are
$t\to \omega(x_{N+k})+t-x_{N+k}$,
$t\to \omega(x_{N+k+1})-t+x_{N+k+1}$).

We now calculate the asymptotics of the probabilities $\mu_n(x_{N+k})$, hoping
to get approximately $(B-A)/n$ (``$\Delta x$'') times the density of
$\mu$ at $x_{N+k}$. To make the argument rigorous, first replace $\mu_n$
by an absolutely continuous version of it, $\mu_n'$, by dispersing the measure
of each $x_{N+k}$ uniformly over the interval $[x_{N+k},x_{N+k+1})$,
and eliminating the measure of the last point $x_{N+n}$.
Since
$\mu_n$ and $\mu_n'$ clearly converge or diverge weakly together
(it will be easy to see from the calculation below that the measure of
the eliminated point $x_{N+n}$ is negligible),
the claim that $\mu_n\to\mu$ thus reduces to checking that the
sequence of  densities $g_n$ of $\mu_n'$ converges to $g(x)$ and
is majorized by an integrable function, so that the dominated convergence
theorem can be applied.

Fix $x\in(A,B)$, and let $k$ such that $x_{N+k}\le x<x_{N+k+1}$
(note that there is an implicit dependence of $n$, of the partition points
$x_{N+j}$ as well as of $k$; as $n$ grows to infinity, $k$ behaves like
$n\cdot(x-A)/(B-A)$). Then:
$$
g_n(x) = \left(\frac{B-A}{n}\right)^{-1}\cdot \mu_n(x_{N+k}) = $$ $$=
  \left(\frac{B-A}{n}\right)^{-1} \cdot
    \prod_{i<N+k} \left(1-\frac{y_i-x_i}{x_{N+k}-x_i}\right) \cdot
    \prod_{i>N+k} \left(1-\frac{x_i-y_{i-1}}{x_i-x_{N+k}}\right) = $$ $$=
\left[ \prod_{i<N} \left(1-\frac{y_i-x_i}{x_{N+k}-x_i}\right) \cdot
       \prod_{i>N+n} \left(1-\frac{x_i-y_{i-1}}{x_i-x_{N+k}}\right) \right]
\cdot$$ $$\cdot \left[ \left(\frac{B-A}{n}\right)^{-1} \cdot
       \prod_{i=N}^{N+k-1} \left(1-\frac{y_i-x_i}{x_{N+k}-x_i}\right) \cdot
       \prod_{i=N+k+1}^{N+n} \left(1-\frac{x_i-y_{i-1}}{x_i-x_{N+k}}\right)
\right]
$$
We treat the two parenthesized expressions 
in the last equation separately: The second one is equal to
$$\left(\frac{B-A}{n}\right)^{-1} \cdot
       \prod_{i=0}^{k-1} \left(1-\frac{y_{N+i}-x_{N+i}}{x_{N+k}-x_{N+i}}
    \right) \cdot
       \prod_{i=k+1}^{n} \left(1-\frac{x_{N+i}-y_{N+i-1}}{x_{N+i}-x_{N+k}}\right)
=$$ $$
\left(\frac{B-A}{n}\right)^{-1} \cdot
\prod_{i=0}^{k-1} \left(1-\frac{(1+d)/2}{k-i}\right) \cdot
\prod_{i=k+1}^n \left(1-\frac{(1-d)/2}{i-k}\right) = $$ $$
\left(\frac{B-A}{n}\right)^{-1} \cdot
\prod_{i=1}^{k} \left(1-\frac{(1+d)/2}{i}\right) \cdot
\prod_{i=1}^{n-k} \left(1-\frac{(1-d)/2}{i}\right) \sim_{n\to\infty}$$ $$
\left(\frac{B-A}{n}\right)^{-1} \cdot
\frac{k^{-(1+d)/2}}{\Gamma((1-d)/2)}\cdot
\frac{(n-k)^{-(1-d)/2}}{\Gamma((1+d)/2)}= $$ $$=
\frac{1}{\pi}\sin(\pi(1+d)/2)\cdot
\left(\frac{k}{n}(B-A)\right)^{-(1+d)/2} 
\cdot \left(\frac{n-k}{n}(B-A)\right)^{-(1-d)/2} = $$ $$
\frac{1}{\pi}\sin(\pi(1+d)/2)\cdot
(x_{N+k}-A)^{-(1+d)/2} 
\cdot (B-x_{N+k})^{-(1-d)/2} \sim $$ $$
\frac{1}{\pi}\sin(\pi(1+d)/2)\cdot
(x-A)^{-(1+d)/2} 
\cdot (B-x)^{-(1-d)/2}
$$ 
It remains to evaluate the asymptotics of the products in the first
parentheses; this is in fact simpler, since in these products the
individual terms tend to zero. In general, a product of the form
$\prod(1+\Delta x_i h(x_i))$ converges to 
$\exp(\int h(u)du)$. In our case, the limit of the two products is then
easily seen to be
$$ \exp\left(
\int_a^A \frac{-(1+\omega'(u))/2}{x-u}du +
\int_B^b \frac{-(1-\omega'(u))/2}{u-x}du \right) $$
Putting the pieces together, we have the formula
$$ \lim_{n\to\infty} g_n(x) =
\frac{1}{\pi}\cos(\pi\omega'(x)/2)\cdot
(x-A)^{-(1+\omega'(x))/2} 
\cdot (B-x)^{-(1-\omega'(x))/2} \cdot $$ $$
\exp\left(
\int_a^A \frac{-(1+\omega'(u))/2}{x-u}du +
\int_B^b \frac{-(1-\omega'(u))/2}{u-x}du \right) $$
We now rearrange the terms slightly, noting that
$$ (x-A)^{-(1+\omega'(x))/2} = (x-a)^{-(1+\omega'(x))/2} \cdot
     \exp\left(\int_a^A \frac{(1+\omega'(x))/2}{x-u}du\right) $$
and
$$ (B-x)^{-(1-\omega'(x)/2)} = (b-x)^{-(1-\omega'(x))/2} \cdot
     \exp\left(\int_B^b \frac{(1-\omega'(x))/2}{u-x}du\right) $$
to finally arrive at
$$ \lim_{n\to\infty} g_n(x) = \frac{1}{\pi} \cos(\pi\omega'(x)/2)\cdot
(x-a)^{-\frac{1+\omega'(x)}{2}} \cdot (b-x)^{-\frac{1-\omega'(x)}{2}} \cdot
$$ $$ \exp\left(\frac{1}{2}
 \int_{[a,A]\cup[B,b]} \frac{\omega'(u)-\omega'(x)}{u-x}du\right)
$$
One more cosmetic change of the formula is to write
$$ \lim_{n\to\infty} g_n(x) = \frac{1}{\pi} \cos(\pi\omega'(x)/2)\cdot
(x-a)^{-\frac{1+\omega'(x)}{2}} \cdot (b-x)^{-\frac{1-\omega'(x)}{2}} \cdot
$$ $$ \exp\left(\frac{1}{2}
 \int_a^b \frac{\omega'(u)-\omega'(x)}{u-x}du\right) = g(x),
$$
since within the segment $[A,B]$ there is no contribution to the integral
inside the exponent.

%We remark at this point that a different way of writing the same expression
%is
%$$ \frac{1}{\pi} \cos(\pi\omega'(x)/2) \cdot [(x-a)(b-x)]^{-1/2} \cdot
% \exp\left(\textrm{P.V.} \int_a^b \frac{\omega'(u)}{u-x}du\right), $$
%where P.V. denotes a principal value integral. It is this form of the formula
%was conjectured by Kerov (see discussion in section 3) as the density
%of the absolutely continuous part of the transition measure, for general
%diagrams.

To complete this part of the proof, we need to show that the densities
$g_n(x)$ are uniformly bounded by some integrable function. Looking back
at the two parenthesized expressions, we see that the first is bounded by
$1$, and the second is
$$
\left(\frac{B-A}{n}\right)^{-1} \cdot
\prod_{i=1}^{k} \left(1-\frac{(1+d)/2}{i}\right) \cdot
\prod_{i=1}^{n-k} \left(1-\frac{(1-d)/2}{i}\right) \le $$ $$ \le
c_1\cdot
\left(\frac{B-A}{n}\right)^{-1} \cdot
\prod_{i=1}^{k+1} \left(1-\frac{(1+d)/2}{i}\right) \cdot
\prod_{i=1}^{n-k} \left(1-\frac{(1-d)/2}{i}\right) \le $$ $$ \le
c_1\cdot
\left(\frac{B-A}{n}\right)^{-1} \cdot
\exp\left(-\sum_{i=1}^{k+1}\frac{(1+d)/2}{i}
          -\sum_{i=1}^{n-k}\frac{(1-d)/2}{i}\right) \le $$ $$ \le
c_2\cdot
\left(\frac{B-A}{n}\right)^{-1} \cdot (k+1)^{-(1+d)/2} \cdot
  (n-k)^{-(1-d)/2} = $$ $$ =
c_2\cdot \left(\frac{(k+1)(B-A)}{n}\right)^{-(1+d)/2} \cdot
         \left(\frac{(n-k)(B-A)}{n}\right)^{-(1-d)/2} \le $$ $$ \le
c_2\cdot (x-A)^{-(1+d)/2} \cdot (B-x)^{-(1-d)/2} $$
for some constants $c_1, c_2$ (depending on $\omega$). Thus we have shown
that 
$$g_n(x) \le c_2 \cdot (x-A)^{-(1+\omega'(x))/2}
  \cdot (B-x)^{-(1-\omega'(x))/2}$$ 
inside any maximal segment of linearity
$[A,B]$, and since $\sup_{x\in[a,b]}|\omega'(x)|<1$, this is an
 integrable function.
Together with the fact that $g_n(x)\to g(x)$ for all $x$ in the interior
of a segment of linearity of $\omega$, this finishes the proof that
$\mu_n\to\mu$ weakly. This proves that $\mu$ is indeed the transition
measure of $\omega$, as was claimed, and therefore (1a) is true for
piecewise linear diagrams.

\subsubsection*{5.1.2. Smooth diagrams}

We now turn to the final approximation step, that of going from piecewise
linear diagrams to piecewise $C^2$ ones. Let $\omega\in{\cal S}[a,b]$,
and define a sequence of approximating piecewise-linear diagrams $\omega_n$,
as follows: for each $n$, partition each of the segments of smoothness of
$\omega$ into $2^n$ equal parts. Then $\omega_n$ is the diagram that is
linear on each of the partition intervals and interpolates $\omega$ at
their endpoints. We denote by ${\cal P}$ the set of all these endpoints
(a countable set), and denote $L= \sup_{x\in[a,b]} |\omega''(x)|<\infty$,
$M = \sup_{x\in[a,b]} |\omega'(x)|<1$.

Let $\mu_n$ be the exterior transition measure, with density $g_n(x)$, 
of $\omega_n$. Let $\mu$ be the measure
whose density $g(x)$ is given by (1a) (for the diagram $\omega$). For
the same reasons as in the previous subsection, it will suffice to prove that
$\mu_n$ converges weakly to $\mu$ as $n\to\infty$, to imply that $\mu$ is
indeed the transition measure of $\omega$. We shall show this in
two steps: first, we show that $g_n(x) \to g(x)$ for \emph{almost all}
$x\in[a,b]$ (somewhat surprisingly, this fails on a large set of
$x$'s, though a set of measure zero). Finally, a suitable boundedness
argument will assure the weak convergence. 

\paragraph{The first step.} Define
$$p(x) = \frac{1}{\pi}\cos(\pi\omega'(x)/2) \cdot
 (x-a)^{-(1+\omega'(x))/2} (b-x)^{-(1-\omega'(x))/2} $$
$$q(x) = \exp\left( \frac{1}{2}
           \int_a^b \frac{\omega'(u)-\omega'(x)}{u-x}du\right) $$
$$p_n(x) = \frac{1}{\pi}\cos(\pi\omega_n'(x)/2) \cdot
 (x-a)^{-(1+\omega_n'(x))/2} (b-x)^{-(1-\omega_n'(x))/2} $$
$$q_n(x) = \exp\left( \frac{1}{2}
           \int_a^b \frac{\omega_n'(u)-\omega_n'(x)}{u-x}du\right), $$
so that $g_n(x)=p_n(x)q_n(x), g(x)=p(x)q(x)$.
Clearly $p_n(x)\to p(x)$ for all $x\in[a,b]\setminus {\cal P}$, we now
try to show $q_n(x) \to q(x)$ (this will fail for some $x$'s, but succeed
for most): For a given $x\in [a,b]\setminus{\cal P}$,
$$ \frac{\omega_n'(u)-\omega_n'(x)}{u-x} \xrightarrow[n\to\infty]{}
   \frac{\omega'(u)-\omega'(x)}{u-x} $$ for all $u\in[a,b]\setminus{\cal P}
 \setminus\{x\}$.
To deduce that $q_n(x)\to q(x)$, some kind of boundedness argument is now
required.
Let $[A,B]$ be the maximal segment of smoothness of $\omega$
containing $x$. Then for all $u\in[a,A]\setminus{\cal P}$ we have
$$ \left|\frac{\omega_n'(u)-\omega_n'(x)}{u-x}\right| \le \frac{2M}{x-A}, $$
and for all $u\in[B,b]\setminus{\cal P}$
$$ \left|\frac{\omega_n'(u)-\omega_n'(x)}{u-x}\right| \le \frac{2M}{B-x},$$
which implies by the dominated convergence theorem that
$$ \int_{[a,A]\cup[B,b]}
     \frac{\omega_n'(u)-\omega_n'(x)}{u-x}du \xrightarrow[n\to\infty]{}
   \int_{[a,A]\cup[B,b]}
     \frac{\omega'(u)-\omega'(x)}{u-x}du$$
Now let $0\le k=k(n) \le 2^n-1$ be such that $x$ is in the $k$th partition
interval of the segment $[A,B]$, i.e. $A+(B-A)k/2^n < x < A+(B-A)(k+1)/2^n$.
We bound $(\omega_n'(u)-\omega_n'(x))/(u-x)$ 
(as a function of $u$)
separately on the different
partition intervals within $[A,B]$. On the interval 
$(A+(B-A)k/2^n,A+(B-A)(k+1)/2^n)$ this expression is zero. On the other
intervals:
We can write $\omega_n'(u) = \omega'(u')$,
$\omega_n'(x) = \omega'(x')$, where $u'$ is in the same partition interval
as $u$ and $x'$ is in the same partition interval as $x$
(since $\omega_n$ is, within each
partition interval, a linear interpolation of $\omega$). Thus if $u$
is in the $j$th interval, then
$$\left|\frac{\omega_n'(u)-\omega_n'(x)}{u-x}\right| =
\left|\frac{\omega'(u')-\omega'(x')}{u'-x'}\right|\cdot
\left|\frac{u'-x'}{u-x}\right| =$$ $$ = 
|\omega''(u'')| \cdot \left| 1+\frac{(u'-u)+(x'-x)}{u-x}\right| \le
L \cdot \left(1+\frac{2(B-A)}{2^n \cdot |u-x|}\right) \le $$ $$ \le
L \cdot \left(1+\frac{2}{|j-k|-1}\right) $$
For $j<k-1$ or $j>k+1$ this gives an effective bound of $3L$. For the
$(k-1)$th and $(k+1)$th interval we are left with the bound of
$(1+2(B-A)/(2^n|u-x|))$, which is not effective at all, since when $x$ and
$u$ are in adjacent intervals they can be arbitrarily close! One may
describe exactly how close they may be using the binary expansion of
$(x-A)/(B-A)$: if we denote by $s_n(x)$ the length of the sequence of
zeroes in this binary expansion starting at the $n$th place, and by
$t_n(x)$ the length of the sequence of ones starting at the $n$th
place, then for $u$ in the $(k\pm 1)$th interval we have 
$|u-x|\ge (B-A) 2^{-(n+s_n(x)\vee t_n(x))}$, and so
$$ \int_{A+(B-A)(k-1)/2^n}^{A+(B-A)k/2^n}
\left|\frac{\omega_n'(u)-\omega_n'(x)}{u-x}\right|du \le $$ $$ \le
\frac{L(B-A)}{2^n}+\frac{2L(B-A)}{2^n}
\int_{A+(B-A)(k-1)/2^n}^{A+(B-A)k/2^n} \frac{du}{x-u} = $$ $$ =
\frac{L(B-A)}{2^n}+\frac{2L(B-A)}{2^n}\log\left(
   \frac{x-(A+(B-A)(k-1)/2^n)}{x-(A+(B-A)k/2^n)}\right) \le $$ $$ \le
\frac{L(B-A)}{2^n}+\frac{2L(B-A)}{2^n}\log\left(\frac{2(B-A)/2^n}{(B-A)/2^{n+s_n(x)\vee t_n(x)}}
 \right)\le $$ $$ \le
\frac{L(B-A)}{2^n}+\frac{4L(B-A)}{2^n}\cdot (s_n(x)\vee t_n(x)).$$
%$$ \int_{A+(B-A)(k-1)/2^n}^{A+(B-A)k/2^n}
%\left|\frac{\omega_n'(u)-\omega_n'(x)}{u-x}\right|du \le $$ $$ \le
%\frac{L(B-A)}{2^n}+\frac{2L(B-A)}{2^n}
%\int_{A+(B-A)(k-1)/2^n}^{A+(B-A)k/2^n} \frac{du}{u-x} = $$ $$ =
%\frac{L(B-A)}{2^n}+\frac{2L(B-A)}{2^n}\log\left(
%   \frac{A+(B-A)k/2^n - x}{A+(B-A)(k-1)/2^n-x}\right) \le $$ $$ \le
%\frac{L(B-A)}{2^n}+\frac{2L(B-A)}{2^n}\log\left(\frac{(B-A)/2^n}{(B-A)/2^{n+s_n(x)\vee t_n(x)}}
% \right)\le $$ $$ \le
%\frac{L(B-A)}{2^n}+\frac{2L(B-A)}{2^n}\cdot (s_n(x)\vee t_n(x)).$$
In a similar manner, integrating on the $(k+1)$th interval gives the
same bound
$$\int_{A+(B-A)(k+1)/2^n}^{A+(B-A)(k+2)/2^n}
\left|\frac{\omega_n'(u)-\omega_n'(x)}{u-x}\right|du \le 
\frac{L(B-A)}{2^n}+\frac{4L(B-A)}{2^n}\cdot(s_n(x)\vee t_n(x)).$$

%, and so our bound becomes
%$$\left|\frac{\omega_n'(u)-\omega_n'(x)}{u-x}\right| \le KL\cdot
%   2^{s_n(x)\vee t_n(x)},$$
%$K$ again being some constant

So, our attempt to prove boundedness of the sequence of integrands failed -
but we are rescued by the fact that it failed on a set of small measure,
namely the two intervals adjacent to the $k$th, and where the values of
the integrands are not too big. In other words, we claim
that, under some further restrictions on $x$, the sequence
$(\omega_n'(u)-\omega_n'(x))/(u-x)$ shall be
 uniformly integrable in $u$. Indeed,
we have proved uniform boundedness on all but the two intervals
adjacent to the $k$th, and on them we have the estimate
$$ \int_{I_{k-1}\cup I_{k+1}} 
 \left|\frac{\omega_n'(u)-\omega_n'(x)}{u-x}\right|du \le
   \frac{2L(B-A)}{2^n}+\frac{8L(B-A)}{2^n}\cdot(s_n(x)\vee t_n(x)),$$
where $I_j = [A+(B-A)j/2^n, A+(B-A)(j+1)/2^n]$ is the $j$th interval.
This bound tends to $0$ (which is what we need to prove uniform integrability)
for those $x\in[A,B]\setminus{\cal P}$ for
which $s_n(x)\vee t_n(x)$ grows asymptotically
at a rate smaller than, say, $2^{n/2}$. But in fact it is a well-known fact
in number theory that
almost every $x$ has this
property (one may prove this directly, or appeal to the stronger
theorem from Feller (1957), p. 197, which says that for almost all 
$z\in[0,1]$, the length $t_n(z)$
of the sequence of zeros in the binary expansion starting at place $n$,
satisfies $\limsup t_n(z)/\log_2(n) = 1$). Thus, for almost every $x$
the sequence $(\omega_n'(u)-\omega_n'(x))/(u-x)$ is uniformly integrable,
and therefore $q_n(x)\to q(x)$, as
was claimed.

\medskip\noindent We note this as a lemma, to be used in section 6:

\paragraph{Lemma 3.} If $f:[a,b]\to\mathbb{R}$ is 
piecewise-continuously-differentiable with bounded derivative,
and $f_n$ is a sequence of piecewise-constant functions obtained by
dividing each interval of differentiability of $f$ into $2^n$ equal
parts and defining $f_n$ on each subinterval as the average value of
$f$ on that subinterval. Then for almost all $x\in[a,b]$
$$
  \int_a^b \frac{f_n(u)-f_n(x)}{u-x}du \xrightarrow[n\to\infty]{}
  \int_a^b \frac{f(u)-f(x)}{u-x}du
$$

\paragraph{The second step.} Having proved $g_n(x)\to g(x)$ for almost
all $x\in[a,b]$, we now finish the proof by showing that 
the $g_n$ are themselves uniformly integrable. Let $x\in[A,B]\subset
 [a,b]$ as before.
The estimates derived above imply that for some constants $k_1, k_2$
(depending on the diagram $\omega$),
$$ g_n(x) \le k_1 (x-A)^{-M} \cdot (B-x)^{-M} e^{k_2 (s_n(x)
\vee t_n(x))/2^n}$$
(If $A\ne a$ and $B\ne b$, then $p_n(x)$ are uniformly bounded by
a constant and $q_n(x)$ is bounded by the above expression; if $A=a$ or
$B=b$, then it is the $p_n(x)$ that contributes the factor
$(x-A)^{-M}$(or, respectively, $(B-x)^{-M}$), which disappears from the
bound on $q_n(x)$.)

To show that this sequence of functions is uniformly integrable on $[A,B]$,
we shall show that it is bounded in L$_p[A,B]$ for some $p>1$. In fact,
$((x-A)(B-x))^{-M}$ is in L$_p[A,B]$ for $1\le p<M^{-1}$, in particular for
$p_0 = (1+M^{-1})/2 > 1$. Let $q_0 = p_0/(p_0-1)$, and let $\epsilon>0$
such that $p_1 = (1+\epsilon)p_0 < M^{-1}$.
If we show that the sequence of functions
$\exp(k_2 (s_n(x)\vee t_n(x))/2^n)$ is uniformly bounded in L$_q$
for any $q\ge 1$,
and thus in particular for $(1+\epsilon)q_0$, then by H\"older's
inequality it will follow that the product of the two expressions,
which majorizes $g_n(x)$, is bounded in L$_{1+\epsilon}$ and
thus uniformly integrable.
And indeed:
$$ \int_A^B \exp(k_2 (s_n(x)\vee t_n(x))/2^n)^q dx = $$ $$ =
 \sum_{j=1}^\infty e^{k_2\cdot q\cdot j/2^n}\cdot
 |\{x\in[A,B]:s_n(x)\vee t_n(x)=j\}| \le $$ $$ \le 2(B-A)\cdot
 \sum_{j=1}^\infty e^{k_2\cdot q\cdot j/2^n}\cdot 2^{-j},$$
and this is finite (and decreasing in $n$, thus bounded) after some
initial value $n=n_0(q)$.

\subsection*{5.2. The interior transition measure}

We now calculate the density of the interior transition measure of $\omega$.
The calculation is quite similar to the one in the previous subsection,
as well as the various proofs of convergence. Therefore, we shall only
write explicitly the calculation of the limiting density for piecewise
linear diagrams. We use the same notation as in section 5.1.1.: $\omega$
is a piecewise linear diagram, $\omega_n$ is the sequence of approximating
rectangular diagrams defined using the equipartition points $x_k$. $\nu_n$
is the interior transition measure of $\omega_n$, and $\nu_n'$ is the
absolutely continuous version of it whereby the probability of each $y_k$
is dispersed uniformly over the interval $[x_k,x_{k+1}]$. $g_n(x)$ is the
density of $\nu_n'$. We calculate: Let $x\in(A,B)\subset[a,b]$,
and let $0\le k< n$ such that $x_{N+k}<x<x_{N+k+1}$, then
$$
g_n(x) = \left(\frac{B-A}{n}\right)^{-1}\cdot \nu(y_{N+k}) = $$ $$ =
\frac{2}{A(\omega_n)}\cdot\left(\frac{B-A}{n}\right)^{-1}
 (x_{N+k+1}-y_{N+k})\cdot (y_{N+k}-x_{N+k})\cdot $$ $$ \cdot 
 \prod_{i<N+k} \left(1+\frac{y_i-x_i}{y_k-y_i}\right)
 \prod_{i>N+k+1}^n \left(1+\frac{x_i-y_{i-1}}{y_{i-1}-y_k}\right) = $$ $$ =
\left[\frac{2}{A(\omega_n)}\cdot
 \prod_{i<N} \left(1+\frac{y_i-x_i}{y_k-y_i}\right) \cdot
       \prod_{i>N+n} \left(1+\frac{x_i-y_{i-1}}{y_{i-1}-y_k}\right)\right]
\cdot $$ $$ \cdot \bigg[
 \left(\frac{B-A}{n}\right)^{-1}
(x_{N+k+1}-y_{N+k})\cdot (y_{N+k}-x_{N+k})\cdot $$ $$ \cdot
\prod_{i=N}^{N+k-1} \left(1+\frac{y_i-x_i}{y_k-y_i}\right)
 \prod_{i=N+k+2}^{N+n} \left(1+\frac{x_i-y_{i-1}}{y_{i-1}-y_k}\right)
\bigg] $$
Again we treat the two parenthesized expressions
 separately; the first one converges
to
$$ \frac{2}{A(\omega)}\cdot
\exp\left(\int_a^A \frac{(1+\omega'(x))/2}{x-u}du +
             \int_B^b \frac{(1-\omega'(x))/2}{u-x}du \right). $$
The second one is
$$
\left(\frac{B-A}{n}\right)^{-1}\cdot \frac{1-\omega'(x)}{2}\cdot
   \frac{B-A}{n} \cdot \frac{1+\omega'(x)}{2}\cdot\frac{B-A}{n} \cdot $$ $$
\cdot \prod_{i=0}^{k-1}\left(1+\frac{(1+\omega'(x))/2}{i}\right) \cdot
 \prod_{i=1}^{n-k-1}\left(1+\frac{(1-\omega'(x))/2}{i}\right) \sim $$ $$
{\Gamma\left(\frac{1+\omega'(x)}{2}\right)}^{-1} \cdot
{\Gamma\left(\frac{1-\omega'(x)}{2}\right)}^{-1} \cdot $$ $$ \cdot
\left(\frac{(B-A)k}{n}\right)^{(1+\omega'(x))/2}\cdot
\left(\frac{(B-A)(n-k)}{n}\right)^{(1-\omega'(x))/2} \sim$$ $$
\frac{1}{\pi}\cos(\pi\omega'(x)/2)\cdot (x-A)^{(1+\omega'(x))/2} \cdot
 (B-x)^{(1-\omega'(x))/2} $$
Altogether we have
$$ \lim_{n\to\infty} g_n(x) =
 \frac{2}{\pi\cdot A(\omega)} \cos(\pi\omega'(x)/2) \cdot
(x-A)^{(1+\omega'(x))/2} \cdot (B-x)^{(1-\omega'(x))/2} \cdot $$ $$ \cdot
\exp\left( \int_a^A \frac{(1+\omega'(x))/2}{x-u}du +
           \int_B^b \frac{(1-\omega'(x))/2}{u-x}du \right) $$
Now as before,
 rearranging the terms and letting $A$ and $B$ tend to $x$ from above
and below gives (1b) for piecewise linear diagrams.

\section*{6. Proof of Theorem 4 and the inversion formulas}

Our main goal in this section is to prove Theorem 3, that gives the
shape of the diagram $\omega$ in terms of the density of the transition
measures of $\omega$. We start by proving Theorem 4, which contains the
essential computational idea behind the proof. We then proceed with the
proof of Theorem 3, where again, the basic idea is to approximate the
diagram by rectangular diagrams, and the transition measures by atomic
measures. There will be two
approximation steps. First, we treat the case of densities which are
step functions with finitely many values. Next we approximate an
arbitrary density satisfying the assumptions of Theorem 3 by such
step functions. The details are given only for the exterior walk case
(Theorem 3a). The proof of Theorem 3b follows the
same reasoning, where the uses of formula (2) and Theorem 5 are replaced
by the their respective analogues, formula (7) and Theorem 6.

\subsection*{6.1. Proof of Theorem 4}

Recall that $p_n(t) = t(t-1)(t-2)...(t-n)$, and 
$(k+\lambda_{n,k})_{k=0}^{n-1}$ are the roots of $p_n'$. Let $0<x<1$, and
denote $k=\lfloor x\cdot n\rfloor$. Then $k+\lambda_{n,k}$ is the root of the
equation
$$ \frac{p_n'(x)}{p_n(x)} = \sum_{j=0}^n \frac{1}{x-j} = 0 $$
in the interval $(k,k+1)$. In other words we have
$$ \sum_{j=0}^k \frac{1}{\lambda_{n,k}+k-j} - 
\sum_{j=k+1}^n \frac{1}{-\lambda_{n,k}+j-k} = 0,$$
or, transforming the indices,
$$ \sum_{j=0}^k \frac{1}{\lambda_{n,k}+j} - 
\sum_{j=0}^{n-k-1}\frac{1}{(1-\lambda_{n,k})+j} = 0.$$
By the classical relations
$$ (13)\qquad \sum_{j=0}^m \frac{1}{u+j} = 
  -\frac{\Gamma'(u)}{\Gamma(u)}+\log(m)+ o(1)_{m\to\infty} $$ $$
(14)\qquad
-\frac{\Gamma'(u)}{\Gamma(u)}+\frac{\Gamma'(1-u)}{\Gamma(1-u)}=
 \pi\cot(\pi u) $$ (equivalent to (11) and (12), respectively)
the latter equation transforms to
$$\pi\cot(\pi \lambda_{n,k}) = \log\left(\frac{n-k-1}{k}\right)+
o(1)_{n\to\infty} = \log\left(\frac{1-x}{x}\right)+o(1)$$ \qed

\subsection*{6.2. The inversion formula for the exterior walk}

\subsubsection*{6.2.1. Step functions}

The notation, and the techniques of approximation,
are much like in the previous sections.
Let $g(x)$, the density of the exterior transition measure $\mu$ of a diagram
$\omega$,
be a step function, taking on finitely many strictly
positive values on $[a,b]$. Thus, $g(x)$ is a mixture of uniform densities
on each of the segments where $g(x)$ is constant. We approximate this
transition measure by the corresponding mixture of discrete uniform measures:
For each $n$, divide as before each (maximal) segment where $g$ is constant
into $n$ equal parts. If $[A,B]$ is one such segment, let
$a=x_0<x_1<x_2<...<x_{l\cdot n}$ be these division points, and let 
$A=x_N < x_{N+1} < ... < x_{N+n} = B$ be the division points
inside the interval $[A,B]$.
(There is an implicit dependence on $n$ here, and as before we suppress it for
convenience of notation.)
Define the measure $\mu_n$ as the discrete measure, concentrated on the
$x_k$, and giving to the point $x_k$ the measure $\mu([x_k,x_{k+1}])$.
Let, for each $n$, $\omega_n$ be the (rectangular) diagram corresponding
to the discrete measure $\mu_n$. Let 
$y_0 < y_1 < y_2 < ... < y_{l\cdot n-1}$ be the sequence of maxima of
the diagram $\omega_n$. For each $k=0,1,2,...,n-1$, since 
$x_{N+k}<y_{N+k}<x_{N+k+1}$, we can write 
$y_{N+k} = x_{N+k}+\lambda_{N+k}\cdot (B-A)/n$ for some $0<\lambda_{N+k}<1$.

It is clear that $\mu_n\to\mu$ weakly as $n\to\infty$, and thus by Theorem 5,
$\omega_n\to\omega$ uniformly on $[a,b]$.
We now proceed to calculate the limit of
$\omega_n$, by calculating the limit of the $\lambda_{N+k}$. By (2),
the $y_i$
are the roots of the equation $\sum_k \frac{\mu_n(x_k)}{x-x_k} = 0$. For
$y_{N+k}$, we write this as
$$ \sum_{j=N}^{N+k} \frac{\mu_n(x_j)}{(N+k-j+\lambda_{N+k})\frac{B-A}{n}} -
  \sum_{j=N+k+1}^{N+n} \frac{\mu_n(x_j)}{(j-(N+k)-\lambda_{N+k})\frac{B-A}{n}}
 = $$ $$ = \sum_{j=1}^{N-1} \frac{\mu_n(x_j)}{x_j-y_{N+k}} +
 \sum_{j=N+n+1}^{l\cdot n} \frac{\mu_n(x_j)}{x_j-y_{N+k}} 
$$
The RHS is $\int_{[a,A]\cup[B,b]} \frac{g(u)du}{u-y_{N+k}} + o_\epsilon(1)$ as
$n\to\infty$, the $o_\epsilon(1)$ being uniformly small for all values of $k$ between
$\epsilon n$ and $(1-\epsilon)n$. The LHS can be rewritten as
$$ g(y_{N+k})\left(
   \sum_{j=0}^k \frac{1}{j+\lambda_{N+k}}-
   \sum_{j=0}^{N-k-1} \frac{1}{j+(1-\lambda_{N+k})} \right).$$
As in the proof of Theorem 4 above, we may use (13) and (14) to transform
this expression as
$$ g(y_{N+k})\left(\pi\cot(\pi \lambda_{N+k})-
      \log\left(\frac{B-y_{N+k}}{y_{N+k}-A}\right)\right) + o_\epsilon(1)$$
(with the same uniformity property). A further rearrangement
of the terms, similar to that done in the previous sections, leads to
the equation
$$ \lambda_{N+k} = \frac{1}{\pi}\textrm{arccot}\left[
   \frac{1}{\pi}\left( \log\left(\frac{b-y_{N+k}}{y_{N+k}-a}\right) + \frac{1}{g(y_{N+k})}
     \int_a^b \frac{g(u)-g(y_{N+k})}{u-x}du\right)\right]+o_\epsilon(1).$$
Now, for $x\in[A,B)$, let $k=k(n)$ such that $x_{N+k} \le x < x_{N+k+1}$, then
$$
\omega_n(x)-\omega(A) = \sum_{j=N}^{N+k-1} (y_j-x_j)-\sum_{j=N+1}^{N+k}
  (x_j-y_{j-1})
 + O\left(\frac{1}{n}\right) = $$ $$
= \sum_{j=N}^{N+k-1} \frac{\lambda_j(B-A)}{n} -
    \sum_{j=N}^{N+k-1}\frac{(1-\lambda_j)(B-A)}{n} +O\left(\frac{1}{n}\right)=
$$ $$ = \frac{2(B-A)}{n} \sum_{j=N}^{N+k-1} \lambda_j - \frac{k(B-A)}{n}+
O\left(\frac{1}{n}\right) = $$ $$ =
 \frac{2(B-A)}{n}
 \left( \sum_{j=N+\epsilon n}^{N+k-1)n} \lambda_j + O(\epsilon) \right)
 -\frac{k(B-A)}{n} + O\left(\frac{1}{n}\right) =  $$ $$ =
 -(x-A) + \frac{2}{\pi}
 \int_A^x \textrm{arccot}\left[
\frac{1}{\pi}\left( \log\left(\frac{b-t}{t-a}\right) + \frac{1}{g(t)}
     \int_a^b \frac{g(u)-g(t)}{u-t}du\right)\right] dt + $$ $$ +
  O\left(\frac{1}{n}\right) + o_\epsilon(1) + O(\epsilon)
$$
Which finishes the proof, since $\epsilon$ was arbitrary.

\subsection*{6.2.2. Piecewise smooth functions}

We now present the final approximation step required to finish the proof
of Theorem 3a.
Let $\omega\in{\cal D}[a,b]$ be such that its exterior transition measure is
absolutely continuous, with a density $g(x)$ that is 
piecewise-continuously-differentiable, has bounded derivative, and is
bounded away from $0$. We approximate $g(x)$ by a sequence $g_n(x)$ of
step functions constructed by the method specified in Lemma 3. Let
$\omega_n\in{\cal D}[a,b]$ be the diagram whose exterior transition measure
is $g_n(x)dx$.
Since $g_n(x)\to g(x)$ for almost all $x\in[a,b]$,
we have $\mu_n \to \mu$ weakly and thus by Theorem 5,
$\omega_n \to \omega$ uniformly on $[a,b]$.
By Lemma 3, and by the fact that the inverse formula holds for the
$g_n$, we have for almost all $t\in[a,b]$
$$ \lim_{n\to\infty} \omega_n'(t)=
-1+\frac{2}{\pi}\textrm{arccot}\left[\frac{1}{\pi}\left(
  \log\left(\frac{b-t}{t-a}\right)
+\frac{1}{g(t)}\int_a^b \frac{g(u)-g(t)}{u-t}du
\right)\right]. $$
(It is here that the boundedness assumptions on $g$ are used.)
Therefore, we have for each $x\in[a,b]$
the chain of equalities
$$ \omega(x) = \lim_{n\to\infty} \omega_n(x) = \lim_{n\to\infty}
 \int_a^x \omega_n'(t)dt = 
 \int_a^x \lim_{n\to\infty} \omega_n'(t)dt$$
So we have shown that
$$ \omega(x) = \int_a^x \left\{ 
-1+\frac{2}{\pi}\textrm{arccot}\left[\frac{1}{\pi}\left(
  \log\left(\frac{b-t}{t-a}\right)
+\frac{1}{g(t)}\int_a^b \frac{g(u)-g(t)}{u-t}du
\right)\right] \right\}dt, $$
and this implies that for almost all $x$
$$ \omega'(x) = 
-1+\frac{2}{\pi}\textrm{arccot}\left[\frac{1}{\pi}\left(
  \log\left(\frac{b-x}{x-a}\right)
+\frac{1}{g(x)}\int_a^b \frac{g(u)-g(x)}{u-x}du
\right)\right],
$$
as was claimed.

%\subsection*{6.3. The inversion formula for the interior walk}
%
%The proof of Theorem 3b follows exactly the same arguments as the proof 
%of Theorem 3a, and we omit it.
%Instead of using formula (2) and Theorem 5, we now use
%formula (7) and Theorem 6.

\section*{7. Other formulas}

We gather in this section some more integral formulas that arise out
of the study of the hook walks. 

\subsection*{7.1. Unrotated diagrams}

We describe the hook walks and formulas (1a) and (1b) for the original
continual diagrams, described simply as increasing functions $f$ on some
interval $[a,b]$ such that $f(a)=0$.
We call such a function an \emph{unrotated diagram}
on $[a,b]$, and denote the set of such diagrams by ${\cal U}[a,b]$.
For each unrotated diagram $f$, there corresponds a diagram
$\omega\in{\cal D}[A,B]$ related to $f$ by
$$ (15) \qquad t = \frac{x+f(x)}{\sqrt{2}},\qquad
\omega(t) = \frac{f(x)+b-x}{\sqrt{2}}, $$
where $A=a$ (for concreteness) and $B=(b+f(b))/\sqrt{2}$.
The domain of $f$ is defined as the set
$$ D_f=\{ (x,y): a\le x\le b, 0\le y\le f(x) \} $$
and the dual domain is
$$ D_f'=\{ (x,y): a\le x\le b, f(x)\le y\le f(b)\}$$
The interior hook of a point $(x,y)\in D_f$ is the set
$$ \{ (x',y')\in D_f: (x'\le x \textrm{\ and\ } y'= y) \textrm{\ or\ }
                      (x'=x\textrm{\ and\ }y'\ge y) \}$$
and the exterior hook of a point $(x,y)\in D_f'$ is
$$ \{ (x',y')\in D_f': (x'=x \textrm{\ and\ } y'\le y) \textrm{\ or\ }
                       (x'\ge x \textrm{\ and\ }y'=y) \} $$
The interior and exterior hook walks are now defined exactly as before.
The interior and exterior transition measures are the distributions of
the $x$-coordinate of the limiting point of the walks.

If $f\in{\cal U}[a,b]$ is continuous,
piecewise twice-continuously-differentiable, with
bounded second derivative and first derivative bounded away from $0$ and
$\infty$, then
the corresponding $\omega$ is in ${\cal S}[A,B]$ and we may use the
change of variables (15) to calculate the density of the transition measures
of $f$. We have:

$$ dt = \frac{1+f'(x)}{\sqrt{2}}dx, \qquad
 \omega'(t) = \frac{d\omega/dx}{dt/dx} = 
\frac{(f'(x)-1)/\sqrt{2}}{(f'(x)+1)/\sqrt{2}} = 1-\frac{2}{1+f'(x)}
$$
%If $g(t)$ is the density of the interior transition measure of $\omega$, 
%and $h(x)$ is the interior transition measure 
%of $f$, then the relation $g(t)dt = h(x)dx$ must
%hold. 
%We calculate:
%$$ g(t)dt = dt\cdot \frac{2}{\pi\cdot A(\omega)} \cdot
%\cos(\pi\omega'(t)/2) \cdot (t-A)^{(1+\omega'(t))/2} \cdot
% (B-t)^{(1-\omega'(t))/2} \cdot $$ $$ \cdot
%\exp\left(-\frac{1}{2}\int_A^B \frac{\omega'(u)-\omega'(t)}{u-t}du\right)
%= $$ $$ =
%dx\cdot \frac{1+f'(x)}{\sqrt{2}} \cdot \frac{2}{\pi A(\omega)}
%\cdot \cos\left(\frac{\pi}{2}\left(1-\frac{2}{f'(x)+1}\right)\right) \cdot
%$$ $$ \cdot
%\left(\frac{x-a+f(x)}{\sqrt{2}}\right)^{\frac{f'(x)}{1+f'(x)}} \cdot
%\left(\frac{b-x+f(b)-f(x)}{\sqrt{2}}\right)^{\frac{1}{1+f'(x)}} \cdot $$ $$
%\cdot \exp\left(\int_a^b 
%  \frac{(1+f'(s))^{-1}-(1+f'(x))^{-1}}{(s-x+f(s)-f(x))/\sqrt{2}} \cdot
%  \frac{1+f'(s)}{\sqrt{2}} ds \right) = $$ $$ =
%dx \cdot \frac{1}{\pi A(f)}(1+f'(x))\cdot 
%\sin\left(\frac{\pi}{1+f'(x)}\right) \cdot $$ $$ \cdot
%(x-a+f(x))^{\frac{f'(x)}{1+f'(x)}}
%\cdot (b-x+f(b)-f(x))^{\frac{1}{1+f'(x)}} \cdot $$ $$ \cdot
%\exp\left(-\int_a^b \frac{1}{s-x+f(s)-f(x)}\cdot
% \frac{f'(s)-f'(x)}{1+f'(x)} ds\right) $$
We leave to the reader to verify:

\paragraph{Theorem 7.} If $f$ satisfies the above conditions, then:

\medskip\noindent {\bf (6a)} The density of the exterior transition
measure for $f$ is equal to
$$ \frac{1}{\pi}(1+f'(x))\cdot
\sin\left(\frac{\pi}{1+f'(x)}\right) \cdot $$ $$ \cdot
(x-a+f(x))^{-\frac{f'(x)}{1+f'(x)}} \cdot
(b-x+f(b)-f(x))^{-\frac{1}{1+f'(x)}} \cdot $$ $$ \cdot
\exp\left(\int_a^b \frac{1}{u-x+f(u)-f(x)} \cdot
  \frac{f'(u)-f'(x)}{1+f'(x)} du \right)
 $$

\medskip\noindent {\bf (6b)} The density of the interior transition
measure for $f$ is equal to
$$ \frac{1}{\pi}(1+f'(x))\cdot
\sin\left(\frac{\pi}{1+f'(x)}\right) \cdot $$ $$ \cdot
(x-a+f(x))^{-\frac{1}{1+f'(x)}} \cdot
(b-x+f(b)-f(x))^{-\frac{f'(x)}{1+f'(x)}} \cdot $$ $$ \cdot
\exp\left(-\int_a^b \frac{1}{u-x+f(u)-f(x)} \cdot
  \frac{f'(u)-f'(x)}{1+f'(x)} du \right)
 $$
In particular, for such $f$ we have
$$ \int_a^b f(x)dx = $$ $$ = \int_a^b \bigg\{
\frac{1}{\pi} (1+f'(x)) \cdot
\sin\left(\frac{\pi}{1+f'(x)}\right) \cdot $$ $$ \cdot
(x-a+f(x))^{\frac{f'(x)}{1+f'(x)}} \cdot
(b-x+f(b)-f(x))^{\frac{1}{1+f'(x)}} \cdot $$ $$ \cdot
\exp\left[ -\int_a^b \frac{1}{u-x+f(u)-f(x)}\cdot
                     \frac{f'(u)-f'(x)}{1+f'(x)} du \right]\bigg\}dx $$
$$ 
\pi = \int_a^b \bigg[
(1+f'(x))\cdot
\sin\left(\frac{\pi}{1+f'(x)}\right) \cdot $$ $$ \cdot
(x-a+f(x))^{-\frac{f'(x)}{1+f'(x)}} \cdot
(b-x+f(b)-f(x))^{-\frac{1}{1+f'(x)}} \cdot $$ $$ \cdot
\exp\left(\int_a^b \frac{1}{u-x+f(u)-f(x)} \cdot
  \frac{f'(u)-f'(x)}{1+f'(x)} du \right)
\bigg] dx
 $$

\subsection*{7.2. The abstract definition of the transition measures}

Equations (3) and (8) may be thought of as an abstract, nonconstructive
way of defining the transition measures of a diagram. Equipped with our
formulas for the densities of the transition measures, we may substitute
them into (3) and (8), respectively, to obtain some more integration
identities:

\paragraph{Theorem 8.} Let $\omega\in {\cal S}[a,b]$. Then for $x\notin[a,b]$:

\medskip\noindent {\bf (7a)}
$$ \exp\left(\frac{1}{2} 
 \int_a^b \frac{\omega'(t)-\textrm{sgn}(t)}{t-x}dt
\right) = $$ $$ =
\int_a^b \bigg[ \frac{1}{\pi} \cos(\pi\omega'(t)/2) \cdot
(t-a)^{-(1+\omega'(t))/2} \cdot (b-t)^{-(1-\omega'(t))/2} \cdot $$ $$ \cdot
\exp\left(\frac{1}{2} \int_a^b \frac{\omega'(u)-\omega'(t)}{u-t}du\right)
\cdot \frac{1}{1-t/x} \bigg] dt
$$

\medskip\noindent {\bf (7b)}
$$ \exp\left(-\frac{1}{2} 
 \int_a^b \frac{\omega'(t)-\textrm{sgn}(t)}{t-x}dt
\right) = $$ $$ =
1-\frac{z(\omega)}{x}-
\int_a^b \bigg[\frac{1}{\pi}\cos(\pi\omega'(t)/2)\cdot
  (t-a)^{(1+\omega'(t))/2} \cdot (b-t)^{(1-\omega'(t))/2} \cdot
$$ $$ \exp\left(-\frac{1}{2}
        \int_a^b \frac{\omega'(u)-\omega'(t)}{u-t} du \right)
\cdot \frac{1}{x(x-t)} \bigg] dt
$$
Note that letting $x\to\infty$ in these equations gives
 (2a) and (2b).

\subsection*{7.3. Relations between walks}

So far, we have only considered the two kinds of hook walks, one of which
leaves from the corner of the dual domain of the diagram, and the other
from a uniformly chosen point in the domain of the diagram. Having calculated
the density for the transition measure of exterior corner walks, 
it is not difficult to transform it into a formula for the density of
the transition measure of an interior
walk leaving from an arbitrary point in the
domain. We do this for unrotated diagrams: Let $f\in{\cal U}[a,b]$ be
continuous,
piecewise twice-continuously-differentiable, with bounded second derivative
and first derivative bounded away from $0$. First, by replacing $f$
with $f^{-1}$ in Theorem 7 we may obtain a formula for the density of
the transition measure of the
\emph{interior corner} walk. Next, a simple scaling transforms this
to a formula for the density of the transition measure for any interior
walk starting from a point $(s,t)$ in the domain. This density, which
we denote by $g_{s,t}(x)$, is given by
$$
(f^{-1}(t)<x<s)\ \ \ 
g_{s,t}(x) = \frac{1}{\pi}(1+f'(x))\sin\left(\frac{\pi}{1+f'(x)}\right)
\cdot $$ $$
\cdot \left(x-f^{-1}(t)+f(x)-t\right)^{-\frac{1}{1+f'(x)}} \cdot
(s-x+f(s)-f(x))^{-\frac{f'(x)}{1+f'(x)}} \cdot $$ $$ \cdot
\exp\left( -\int_{f^{-1}(t)}^s \frac{1}{u-x+f(u)-f(x)}\cdot
  \frac{f'(u)-f'(x)}{1+f'(x)} du\right)
$$
We can now write some equations that describe some of the interrelations
between the different walks: The first equation expresses the defining
fact that each step of the walk goes from the current point to a point
in the hook of the current point, chosen uniformly. Thus, the densities
$g_{s,t}$ must satisfy
\paragraph{Theorem 9.}
$$ g_{s,t}(x) = \frac{1}{s-f^{-1}(t)+f(s)-t}
 \left(\int_x^s g_{v,t}(x)dv + \int_t^{f(x)} g_{s,v}(x)dv\right)$$
(This equation is equal in content, but not in form, to eq. (4.3.5) of 
Kerov (1993).)

\medskip\noindent
The second equation expresses the fact that the uniform interior walk is really
a mixture of all the walks with different given starting points $(s,t)$ ,
with respect to the normalized area measure $dsdt/A(f)$.
This implies the identity
\paragraph{Theorem 10.}
$$
\frac{1}{\pi A(f)}(1+f'(x)) \sin\left(\frac{\pi}{1+f'(x)}\right) \cdot
$$ $$ \cdot
(x-a+f(x))^{\frac{f'(x)}{1+f'(x)}}\cdot (b-x+f(b)-f(x))^{\frac{1}{1+f'(x)}}
\cdot $$ $$ \cdot
\exp\left(-\int_a^b\frac{1}{u-x+f(u)-f(x)}\cdot\frac{f'(u)-f'(x)}{u-x}du
\right) = $$ $$ =
\int_x^b \int_0^{f(x)} g_{s,t}(x) \frac{dt\ ds}{A(f)} $$
Which, after cancelling identical terms on both sides, becomes
$$
(x-a+f(x))^{\frac{f'(x)}{1+f'(x)}} \cdot
(b-x+f(b)-f(x))^{\frac{1}{1+f'(x)}} \cdot $$ $$ \cdot
\exp\left(-\int_a^b\frac{1}{u-x+f(u)-f(x)}\cdot\frac{f'(u)-f'(x)}{u-x}du
\right) = $$ $$
\int_x^b \int_0^{f(x)} \bigg[ 
\left(x-f^{-1}(t)+f(x)-t\right)^{-\frac{1}{1+f'(x)}} \cdot
(s-x+f(s)-f(x))^{-\frac{f'(x)}{1+f'(x)}} \cdot $$ $$ \cdot
\exp\left(-\int_{f^{-1}(t)}^s \frac{1}{u-x+f(u)-f(x)}\cdot
  \frac{f'(u)-f'(x)}{u-x} du\right) \bigg] dt\ ds
$$

\section*{References}

Cifarelli, D.M., Regazzini, E. (1990). Distribution functions of means of a
Dirichlet process. \textit{Ann. Statist.} $\mathbf{18}$, no. 1,
429-442.

\medskip\noindent
Diaconis, P., Kemperman, J. (1996). Some new tools for Dirichlet priors.
In: Bayesian Statistics 5, ed. J.M. Bernardo, J.O. Berger, A.P. Dawid,
F.M. Smith. Oxford University Press.

\medskip\noindent
Feller, W. (1957). An Introduction to Probability Theory and its
Applications, vol. 1, 2nd ed. Wiley, New York.

\medskip\noindent
Greene, C., Nijenhuis, A., Wilf, H. (1979). A probabilistic proof of a formula
for the number of Young tableaux of a given shape. \textit{Adv. Math.}
$\mathbf{31}$, no. 1, 104-109.

\medskip \noindent
Greene, C., Nijenhuis, A., Wilf, H. (1984). Another probabilistic proof in
the theory of Young tableaux. J. Combin. Theory Ser. A $\mathbf{37}$,
no. 2, 127-135.

\medskip\noindent
Ivanov, V., Olshanski, G. (2001). Kerov's central limit theorem for
the Plancherel measure on Young diagrams. In: Symmetric Functions
2001: Surveys of Developments and Perspectives. Proc. NATO Advanced
Study Institute, ed. S. Fomin, Kluwer, 2002.

\medskip \noindent
Kerov, S.V. (1993). Transition probabilities for continual Young diagrams
and the Markov moment problem. 
(Russian) \textit{Funktsional. Anal. i Prilozhen.}
 $\mathbf{27}$, no. 2, 32-49, 96;
translation in \textit{Funct. Anal. Appl.} $\bf{27}$, no. 2, 104-117.

\medskip\noindent
Kerov, S.V. (1998). Interlacing measures. In: Kirillov's seminar on
representation theory, 35-83, \textit{Amer. Math. Soc. Transl. Ser. 2}, 181,
Amer. Math. Soc., Providence, RI.

\medskip \noindent
Kerov, S.V. (1999). A differential model of growth of Young diagrams.
(Russian) Proceedings of the St. Petersburg Mathematical Society, vol. IV,
111-130; translation in \textit{Amer. Math. Soc. Transl. Ser. 2}, 188,
Amer. Math. Soc., Providence, RI.

\medskip \noindent
Krein, M.G., Nudelman, A.A. (1977)., The Markov moment problem and extremal
problems, AMS, Providence, RI.

\medskip \noindent
Pittel, B. (1986). On growing a random Young tableau. 
\textit{J. Comb. Th. Ser. A} $\bf{41}$, no. 2., 278-285.

\bigskip\noindent
Author's current address:\\
\textsc{Department of Mathematics}\\
\textsc{Weizmann Institute of Science}\\
\textsc{Rehovot 76100}\\
\textsc{ISRAEL}\\

\medskip\noindent e-mail: \texttt{romik@wisdom.weizmann.ac.il}

\end{document}